\documentclass[12pt,reqno]{amsart}
\usepackage{amsmath}
\usepackage{amssymb}

\author[Johan H. Tykesson]{Johan H. Tykesson}
\address[Johan H. Tykesson]
{Department of Mathematics, Division of Mathematical Statistics,
Chalmers University of Technology and G\"{o}teborg University, S-412
96 } \email[Johan H. Tykesson]{johant@math.chalmers.se}

\usepackage{latexsym}
\usepackage{amsfonts}
\usepackage{amsmath}
\usepackage{amsthm}
\usepackage{epsfig}

\numberwithin{equation}{section} \numberwithin{figure}{section}
\newtheorem{thm}{Theorem}

\numberwithin{thm}{section}
\newtheorem{lemma}[thm]{Lemma}
\newtheorem{cor}[thm]{Corollary}
\newtheorem{definition}[thm]{Definition}

\def \R {\mathbb R}

\def \hyp {{\mathbb H}^2}
\def \hypn {{\mathbb H}^n}
\def \P {{\bf P}}
\def \E {{\bf E}}

\begin{document}

\title[Uniqueness in the Poisson Boolean model] {Continuum percolation
  at and above the uniqueness treshold on homogeneous spaces}

\begin{abstract}
  We consider the Poisson Boolean model of continuum percolation on a
  homogeneous space $M$. Let $\lambda$ be the intensity of the
  underlying Poisson process. Let $\lambda_u$ be the infimum of the
  set of intensities that a.s. produce a unique unbounded component.
  First we show that if $\lambda>\lambda_u$ then there is a.s. a
  unique unbounded component at $\lambda$. Then we let $M=\hyp\times
  \R$ and show that at $\lambda_u$ there is a.s. not a unique
  unbounded component. These results are continuum analogies of
  theorems by H\"aggstr\"om, Peres and Schonmann.
\end{abstract}


\maketitle \noindent



\section{Introduction and results}

In this paper we show continuum analogies to some theorems
concerning the uniqueness phase in the theory of independent bond
and site percolation on graphs. Before turning to our results, we review these
theorems.

Let $G=(V,E)$ be an infinite transitive graph with vertex set $V$ and
edge set $E$. Keep each edge with probability $p$ and delete it
otherwise, independently for all edges. We call this independent bond
percolation on $G$ at level $p$, and let $\P_p$ be the corresponding
probability measure on the subgraphs of $G$. A connected component in
the random subgraph obtained in percolation is called a cluster. Let
\[p_c(G):=\inf\{p\,:\,\P_p-\mbox{a.s. there is an infinite
  cluster}\}\] be the critical probability for percolation.

In what follows we will discuss percolation at different levels, and
when we do this, we always use the following coupling. To each $e\in
E$ we associate an independent random variable $U_e$ which is
uniformly distributed on $[0,1]$. Then say that $e$ is kept at level
$p$ if $U_e<p$ and deleted otherwise. Using this construction, we
have that if $p_1<p_2$ then any edge kept at level $p_1$ is also
kept at level $p_2$. Therefore we call this coupling the monotone
coupling.

Now suppose that $p_c<p_1<p_2$ and use the monotone coupling. We say
that an infinite cluster at level $p_2$ is $p_1$-stable if it
contains an infinite cluster at level $p_1$.  H\"aggstr\"om and Peres
\cite{haggperes} showed the following theorem:

\begin{thm}\label{hpsats} Suppose $G$ is a transitive unimodular graph and that
  $p_c(G)<p_1<p_2\le 1$. Then any infinite cluster at level $p_2$ is a.s.
  $p_1$-stable.
\end{thm}

The proof of \ref{hpsats} relies on a technique called the
mass transport principle, which is not available in the
non-unimodular setting. However, Schonmann \cite{schonmann} was able
to avoid the use of the mass transport principle and showed:

\begin{thm}\label{schonmannsats} Suppose $G$ is a transitive graph and that
  $p_c(G)<p_1<p_2\le 1$. Then any infinite cluster at level $p_2$ is a.s.
  $p_1$-stable.
\end{thm}

Theorem \ref{schonmannsats} has the following immediate consequence.
Let \[p_u(G):=\inf\{p\,:\,\P_p-\mbox{a.s. there is a unique infinite
  cluster}\}\] be the uniqueness treshold for percolation.

\begin{cor}\label{schonmanncor} Suppose $G$ is a transitive graph and
  that $p>p_u(G)$. Then \newline $\P_p[\mbox{there is a unique
    infinite cluster}]=1$.\end{cor}

So Corollary \ref{schonmanncor} settles what happens above $p_u$.
But there is also the question what happens at $p_u$. It turns out
that the answer depends on the graph. The following theorem of Peres
\cite{peres} is of special interest to us:

\begin{thm}\label{peressats} Let $G=(V_G,E_G)$ and $H=(V_H,E_H)$ be
  two infinite transitive graphs and suppose $G$ is
  nonamenable and unimodular. Then at $p_u(G\times H)$ there is a.s.
  not a unique unbounded component.\end{thm}

In contrast to this result, Benjamini and Schramm \cite{itai1}
showed that on any planar, transitive unimodular graph with one end,
there is a.s. a unique infinite cluster at $p_u$.

We will now discuss analogues of Theorems \ref{hpsats},
\ref{schonmannsats} and \ref{peressats} in a continuum percolation
setting. A Riemannian manifold $M$ is said to be a (Riemannian)
homogeneous space if for each $x,y\in M$ there is an isometry that
takes $x$ to $y$. Throughout this paper we assume that $M$ is an
unbounded homogeneous space, with metric $d_M$ and volume measure
$\mu_M$.  When it is clear which space we are working with we will
write $d=d_M$ and $\mu=\mu_M$. We let $0$ denote the origin of the
space.

For one of the main results below it is possible to give a shorter
proof under the additional assumption that $M$ is a symmetric space. A
connected Riemannian manifold $M$ is said to be a (Riemannian)
symmetric space if for each point $p\in M$ there is an isometry $I_p$
such that $I_p(p)=p$ and $I_p$ reverses geodesics through $p$.  The
most important symmetric spaces where it makes sense to study
continuum percolation are arguably $n$-dimensional Euclidean space
${\mathbb R}^n$ and $n$-dimensional hyperbolic space $\hypn$.  Also
products of symmetric spaces are symmetric spaces, for example
$\hyp\times {\mathbb R}$. Any symmetric space is homogeneous. For an
example of a noncompact space which is homogeneous but not symmetric,
one may consider certain Damek-Ricci spaces, see \cite{berndt}. Next
we introduce the Poisson Boolean model of continuum percolation.

Let $S(x,r):=\{y\in M\,:\,d_M(x,y)\le r\}$ be the closed ball with
radius $r$ centered at $x$. Let $X^{\lambda}$ be a Poisson point process on
$M$ with intensity $\lambda$.  Around every point of $X^{\lambda}$ we
place a ball of unit radius, and denote by $C^{\lambda}$ the region of
the space that is covered by some ball, that is
$C^{\lambda}:=\cup_{x\in X^{\lambda}}S(x,1)$. We remark that all
proofs below work if we instead consider the model with some arbitrary fixed
radius $R$.  Write $\P_{\lambda}$ for the probability measure
corresponding to this model, which is called the Poisson Boolean model
with intensity $\lambda$.

Next we introduce some additional notation. Let
$V^{\lambda}:=(C^{\lambda})^c$ be the vacant region. Let
$C^{\lambda}(x)$ be the component of $C^{\lambda}$ containing $x$.
$C^{\lambda}(x)$ is defined to be the empty set if $x$ is not covered.
Let $X^{\lambda}(A)$ be the Poisson points in the set $A$.
Furthermore denote by $C^{\lambda}[A]$ the union of all balls centered
within the set $A$.  With $N_C$ and $N_V$ we denote the number of
unbounded connected components of $C^{\lambda}$ and $V^{\lambda}$
respectively. The number of unbounded components for the Poisson
Boolean model on a homogeneous space is an a.s. constant which equals
$0$, $1$ or $\infty$. The proof of this is very similar to the
discrete case, see for example Lemma 2.6 in \cite{jonasson2}. As in
the discrete case, we introduce two critical intensities. Let
\[\lambda_c(M):=\inf\{\lambda\,:\,N_C>0\mbox{ a.s.}\} \mbox{ and
}\lambda_u(M):=\inf\{\lambda\,:\,N_C=1\mbox{ a.s.}\}\]

be the critical intensity for percolation and the uniqueness treshold
for the Poisson Boolean model. \vspace{1em}

{\bf Remark.} Obviously it is only interesting to study what happens
at and above $\lambda_u$ when $\lambda_u<\infty$. For example this is
case for $\hyp\times {\mathbb R}$ and may be proved by adjusting the
arguments for the $\hyp$ case, see \cite{tykesson}. Simple
modifications (just embed a different graph in the space) of the
arguments in Lemma 4.8 in \cite{tykesson} shows that for $\lambda$
large enough there are a.s.  unbounded components in $C^{\lambda}$ but
a.s. no unbounded components in $V^{\lambda}$.  Since any two
unbounded components in $C^{\lambda}$ must be separated by some
unbounded component in $V^{\lambda}$ it follows that for $\lambda$
large enough there is a.s. a unique unbounded component in
$C^{\lambda}$.\vspace{1em}

We will often work with the model at several different intensities at
the same time. Suppose we do this at the intensities
$\lambda_1<\lambda_2<...<\lambda_n$. Then we will always assume that
$C^{\lambda_{i+1}}$ is the union of $C^{\lambda_i}$ and balls centered
at the points of a Poisson process with intensity
$\lambda_{i+1}-\lambda_{i}$. We call this the monotone coupling and is
obviously the analogy of the discrete coupling described earlier.

Now suppose $\lambda_1<\lambda_2$ and use the monotone coupling. We
say that an unbounded component in $C^{\lambda_2}$ is
$\lambda_1$-stable if it contains some unbounded component in
$C^{\lambda_1}$. We now state a continuum version of Theorem
\ref{hpsats}.

\begin{thm}\label{uniqmon1}
  Consider the Poisson Boolean model on the homogeneous space $M$. Suppose
  $\lambda_c(M)<\lambda_1<\lambda_2<\infty$. Then a.s. any unbounded
  $\lambda_2$-component is $\lambda_1$-stable.
\end{thm}

From Theorem \ref{uniqmon1}, the following corollary is immediate.

\begin{cor}\label{uniqcor1}
Consider the Poisson Boolean model on the homogeneous space M. Suppose
$\lambda_u(M)<\lambda$. Then $\P_{\lambda}[N_C=1]=1$.
\end{cor}

\noindent {\bf Remark.} Corollary \ref{uniqcor1} is known in the cases
$M={\mathbb R}^n$ for any $n\ge 2$ (see \cite{meester}) and $M=\hyp$
(see \cite{tykesson}).

We will present two proofs of Theorem \ref{uniqmon1}. The first is
inspired by the proof of Theorem \ref{hpsats} and the second is
inspired by the proof of Theorem \ref{schonmannsats}. To get a
continuum analogy to Theorem \ref{peressats} we consider the Poisson
Boolean model on a product space.

\begin{thm}\label{htimesrsats} Consider the Poisson-Boolean model on
  $\hyp\times \R$. At $\lambda_u$ there is a.s. not a unique
  unbounded component.
\end{thm}

Note that if one instead considers the model on $\hyp$, then
Corollary 5.10 in \cite{tykesson} says that at $\lambda_u$ there is
a.s. a unique unbounded component. We now move on to the proofs.

\section{Uniqueness monotonicity}

In this section we first present a short proof for Theorem
\ref{uniqmon1} in the symmetric case, and then a proof which only
needs the assumption that the space is homogeneous.

First we present an essential ingredient to the first proof, the mass
transport principle which is due to Benjamini and Schramm
\cite{itai1}. We denote the group of isometries on the symmetric space $M$ by
Isom($M$).

\begin{definition}
\label{diagonalinv} A measure $\nu$ on $M\times M$ is said to be
\emph{diagonally
  invariant} if for all measurable $A,\,B\subset M$ and
$g\in$\emph{Isom(M)} \[\nu(g A\times g B)=\nu(A\times
B).\]\end{definition}

\begin{thm}({\sc Mass Transport Principle on M})
\label{masstransport} If $\nu$ is a positive diagonally invariant
measure on $ M\times M$ such that $\nu(A\times M)<\infty$ for some
open $A\subset M$, then
\[\nu(B\times M)=\nu(M \times B)\] for all measurable $B\subset
M$.
\end{thm}

Actually the mass transport principle is proved in \cite{itai1} for
the case when $M=\hyp$, but as is remarked there, it holds for any
symmetric space.\vspace{1em}

{\em Proof of Theorem \ref{uniqmon1} in the symmetric case:} Suppose
$\lambda_c<\lambda_1<\lambda_2$. We couple $C^{\lambda_1}$ and
$C^{\lambda_2}$ using the monotone coupling. We are done if we can
show that any unbounded component of $C^{\lambda_2}$ contains an
unbounded component of $C^{\lambda_1}$. Since any ball in
$C^{\lambda_1}$ is also present in $C^{\lambda_2}$, this is
equivalent to show that any unbounded component of $C^{\lambda_2}$
intersects an unbounded component of $C^{\lambda_1}$. For any point
$x\in M$ let
  \[D(x):=\inf\{d(x,y)\mbox{ : }y\mbox{ is in an unbounded component
    of }C^{\lambda_1}\}\] and let
\[\tilde{D}(x):=\left\{\begin{array}{ll} \inf_{y\in
    C^{\lambda_2}(x)}D(y), & \mbox{if }x\in C^{\lambda_2} \\ D(x), \mbox{ otherwise}
\end{array} \right.\]
Define the random set $H$ to be the set of all points $x$
satisfying the conditions
\begin{itemize}
\item $C^{\lambda_2}(x)$ is a $\lambda_1$-unstable unbounded
  component
\item $D(x)\le \tilde{D}(x)+1/2$
\end{itemize}
and write $B(x)$ for the event that $x\in H$. Suppose that
$C^{\lambda_2}$ contains an unbounded component which does not
intersect an unbounded component of $C^{\lambda_1}$. Then this
unbounded component contains regions of positive volume in
$H$, so it suffices to show that $\P[B(x)]=0$. Let $H(x)$ be the
connected component of $H$ containing $x$. Let
$B^{\infty}(x):=B(x)\cap \{\mu(H(x))=\infty\}$ and $B^f
(x):=B(x)\cap\{\mu(H(x))<\infty\}$. The events $B^f$ and
$B^{\infty}$ partition $B$. First we show that $\P[B^f(x)]=0$ using
the mass transport principle.

In any unbounded component of $C^{\lambda_2}$ not intersecting an
unbounded component of $C^{\lambda_1}$ we put mass of unit density.
Then all mass in the unbounded component is transported to the
regions in the unbounded component which are in $H$. Let
$\nu(A\times B)$ be the expected mass sent from the set $A$ to the
set $B$. Then $\nu$ is easily seen to be a positive diagonally
invariant measure on $M\times M$. If $\P[B^f(x)]>0$ then if $A$ is
some connected set of finite positive volume, $A$ will get an
infinite amount of incoming mass with positive probability, that is
$\nu(M\times A)=\infty$. On the other hand, $\nu(A\times M)$, the
amount of mass going out from $A$, is at most $\mu(A)<\infty$. Thus
by the mass transport principle $\P[B^f(x)]=0$.

Next we show $\P[B^{\infty}(x)]=0$ by showing
$\P[B^{\infty}(x)|\tilde{D}(x)=r]=0$ for any $r$. Fix $r$. Suppose
$\{\tilde{D}(x)=r\}$ happens. Then for $B^{\infty}(x)$ to happen,
there must be infinitely many balls in $C^{\lambda_2}(x)$ centered at
distance between $r+1$ and $r+1+1/2$ from unbounded components in
$C^{\lambda_1}$. However, this is not possible, as is seen by
``building'' up the process as follows. Condition on $C^{\lambda_1}$
and then on those balls in $C^{\lambda_2}$ that are centered at
distance at least $r+1$ from unbounded components in $C^{\lambda_1}$.
We have then not conditioned on the balls that are not present in
$C^{\lambda_1}$ but in $C^{\lambda_2}$, and centered at a distance
between $0$ and $r+1$ from unbounded components of $C^{\lambda_1}$.
These balls are centered at a Poisson process of intensity
$\lambda_2-\lambda_1>0$ in this region, and this Poisson process is
independent of everything else we have previously conditioned on.
Thus if there are infinitely many balls in $C^{\lambda_2}(x)$ centered
at distance between $r+1$ and $r+3/2$ from unbounded components in
$C^{\lambda_1}$, then balls centered at the points of the previously
mentioned Poisson process will almost surely connect
$C^{\lambda_2}(x)$ to some unbounded component in $C^{\lambda_1}$.
Thus $\P[B^{\infty}(x)|\tilde{D}(x)=r]=0$ for any $r$ and consequently
$\P[B^{\infty}(x)]=0$. {$\Box$ \vspace{1em}}

For the second proof of Theorem \ref{uniqmon1}, we need some
preliminary results. First we describe a method to find the
component of $C^{\lambda}$ containing $x$. This may be considered to
be the continuum version of the algorithm described in for example
\cite{schonmann} for finding the cluster of a given vertex in
discrete percolation.

At $x$, we grow a ball with unit speed until it has radius $1$, when
the growth of the ball stops. Whenever the boundary of this ball
hits a Poisson point, a new ball starts to grow with unit speed at
this point until it has radius $2$. In the same way, every time a
new Poisson point (which has not already been found) is hit by the
boundary of a growing ball, a ball starts to grow at this point
until it has radius $2$ and so on. Let $L_t^{\lambda}(x)$ denote the
set which has been passed by the boundary of some ball at time $t$.
If $C^{\lambda}(x)$ is bounded, then $L_t^{\lambda}(x)$ stops
growing at some random time $T$. In this case
$C^{\lambda}[L_T^{\lambda}(x)]=C^{\lambda}(x)$ and
$L_T^{\lambda}(x)$ is the $1$-neighbourhood of $C^{\lambda}(x)$. (If
the first ball does not hit any Poisson point, then $C^{\lambda}(x)$
is the empty set). If $C^{\lambda}(x)$ is unbounded, then
$L_t^{\lambda}(x)$ never stops growing. We will refer to this
procedure to as "growing the component containing $x$".

In what follows we will make use of the following lemma, which may
be considered intuitively clear. The proof is inspired by the proof
of the corresponding lemma for the discrete situation which is Lemma
1.1 of \cite{schonmann}.

\begin{lemma}\label{ballemma}
  Consider the Poisson Boolean model on a homogeneous space $M$. Let $R>0$
  and let $\lambda>\lambda_c$. Any unbounded component of
  $C^{\lambda}$ contains balls of radius $R$.
\end{lemma}

\noindent For the proof we need to introduce some further notation.
For a connected set $A$ containing $x$ we let $C^{\lambda}(x,A)$ be
all points in $A$ which can be connected to $x$ by some curve in
$C^{\lambda}\cap A$. Let $E_r(x)$ be the union of all balls centered
within $S(x,r+1)$ that are connected to $x$ via a chain of balls
centered within $S(x,r+1)$. Note that $C^{\lambda}(x,S(x,r))\subset
E_r(x)$.


Let $\delta_r(x):=\sup_{y\in E_r(x)\setminus S(x,r)}d(y,\partial
S(x,r))$ where the supremum is defined to be $0$ if $E_r(x)\setminus
S(x,r)$ is the empty set.  Let $\{A\leftrightarrow B\}$ be the event
that there is some continuous curve in $C^{\lambda}$ which intersects
both the set $A$ and the set $B$. Let $A^o$ be the interior of the set
$A$.

\begin{proof}

  Fix a point $x\in M$. Since the case $R\le 1$ is trivial, we suppose
  $R>1$. For any $r>0$ let $F_r(x):=\{x\leftrightarrow
  \partial S(x,r)\}$ and let \[G_r(x):=\{C^{\lambda}(x,S(x,r))\mbox{
    does not contain a ball of radius }R\}.\] Let $D_r(x):=F_r(x)\cap
  G_r(x)$. Let $D(x)$ be the event that $x$ is an unbounded component
  that does not contain a ball of radius $R$.  Then $D_r(x)\downarrow
  D(x)$ so it is enough to show that $\P[D_r(x)]\rightarrow 0$ as
  $r\rightarrow \infty$. Note that $D_r(x)$ is independent of the
  Poisson process outside $S(x,r+1)$.  Also note that
  $\delta_{r}(x)\in[0,2]$.

  If $D_r(x)\cap \{\delta_{r}(x)< 1/2\}$ occurs, then there is a ball
  centered in $S(x,r-1/2)^o\setminus S(x,r-1)^o$ which is connected to
  $x$ by a chain of balls centered in $S(x,r-1/2)^o$. All these
  balls are also included in the set $E_{r-1/2}(x)$, and one of these
  balls is centered at a distance at most $1/2$ from $\partial
  S(x,r-1/2)$. This gives
  \begin{equation}\label{incleq}D_r(x)\cap \{\delta_r(x)<1/2\}\subset
    D_{r-1/2}(x)\cap \{\delta_{r-1/2}(x)\ge 1/2\}.\end{equation}
 
  We will now proceed by contradiction. Suppose that $\P[D(x)]>0$ and
  that $\lim_{r\rightarrow \infty}\P[\delta_{r}(x)<1/2|D_r(x)]=1.$
  These assumptions imply that \begin{multline*}\lim_{r\rightarrow
      \infty}\P[D_r(x)\cap
    \{\delta_r(x)<1/2\}]\\=\lim_{r\rightarrow\infty}\P[\delta_{r}(x)<1/2|D_r(x)]\P[D_r(x)]=\lim_{r\rightarrow\infty}\P[D_r(x)]=\P[D(x)]>0.\end{multline*}
  However, by (\ref{incleq}) we get that
  \begin{multline*}\limsup_{r\rightarrow
      \infty}\P[\delta_{r-1/2}(x)\ge
    1/2|D_{r-1/2}(x)]\ge\limsup_{r\rightarrow
      \infty}\P[D_{r-1/2}(x)\cap\{\delta_{r-1/2}(x)\ge
    1/2\}]\\\ge\lim_{r\rightarrow\infty}\P[D_r(x)\cap
    \{\delta_r(x)<1/2\}]>0,\end{multline*} so that in particular
  $\P[\delta_r(x)\ge 1/2|D_r(x)]$ does not go to $0$ as $r\rightarrow
  \infty$ which contradicts the assumption $\lim_{r\rightarrow
    \infty}\P[\delta_{r}(x)<1/2|D_r(x)]=1.$ Thus we conclude that
  $\P[D(x)]=0$ or/and $\liminf_{r\rightarrow
    \infty}\P[\delta_r(x)<1/2|D_r(x)]<1.$ We now assume
  $\liminf_{r\rightarrow \infty}\P[\delta_r(x)<1/2|D_r(x)]<1$ and show
  that this implies $\P[D(x)]=0$. By the assumption, we may pick a
  constant $c_1>0$ and a sequence of positive numbers
  $\{a_k\}_{k=1}^{\infty}$ such that $a_{k+1}-a_k\ge 2R+1$ and
  $\P[\delta_{a_k}(x)\ge 1/2|D_{a_k}(x)]\ge c_1$ for all $k$. 
  On the event $D_{a_k}(x)$ we may pick a point $Y$ on $\partial
  S(x,a_k+R+1)$ such that if $S(Y,R+\max(0,1-\delta_{a_k}(x)))$ is
  completely covered by balls centered within $S(Y,R)$, then
  $D_{a_{k+1}}(x)^c$ occurs since a ball of radius $R$ has been found
  in $C(x,S(x,a_{k+1}))$ (this ball is contained in
  $C(x,S(x,a_{k+1}))$ since $a_{k+1}-a_k\ge 2R+1$ and $R>1$). The
  configuration of balls within $S(Y,R)$ is independent of the Poisson
  process within $S(x,a_k+1)$.  Now let $\Delta_k$ be a random
  variable with the same distribution as the conditional distribution
  of $\delta_k(x)$ given the event $D_k(x)$.  By the above
  observations we get that
  \[\P[D_{a_{k+1}}(x)^c|D_{a_k}(x)]\ge\P[S(0,R+\max(0,1-\Delta_k))\subset
  C^{\lambda}[S(0,R)]]\ge c_2\] for some constant $c_2>0$ for all $k$.
  This implies $\lim_{k\rightarrow \infty}\P[D_{a_k}(x)]=0$ and
  consequently $\P[D(x)]=0$.\end{proof}

{\sl Proof of Theorem \ref{uniqmon1}:}

We consider the monotone coupling of the model at intensities
$\lambda_1<\lambda_2$, and we write
$C=(C^{\lambda_1},C^{\lambda_2})$. Let
\[E(x):=\{x\mbox{ is in an unbounded $C^{\lambda_2}$ component which
  is $\lambda_1$-unstable}.\}\]
\noindent Let \[E_1(x):=E(x)\cap \{\tilde{D}(x)\le 3\}\mbox{ and
}E_2(x):=E(x)\cap\{\tilde{D}(x)>2\},\] where $\tilde{D}$ is defined as
in the proof of Theorem \ref{uniqmon1}.

Finally let $E$, $E_1$ and $E_2$ be the events that $E(x)$, $E_1(x)$
and $E_2(x)$ respectively happen for some $x$.

We will first show that $\P[E_2(x)]=0$. Pick $a$ and $R=R(a)$ so
that
\[\P[S(x,R)\mbox{ intersects an unbounded component of
  $C^{\lambda_1}$}]\ge 1-a\]

Let $Z^{'}=(Z^{'\lambda_1},Z^{'\lambda_2})$ and
$Z^{''}=(Z^{''\lambda_1},Z^{''\lambda_2})$ be two independent copies
of $C$, and let $X^{'}=(X^{'\lambda_1},X^{'\lambda_2})$ and
$X^{''}=(X^{''\lambda_1},X^{''\lambda_2})$ be their underlying
Poisson processes. A prime will be used to denote objects relating
to $Z^{'}$ and a double prime will be used to denote objects
relating to $Z^{''}$.

Grow the component of $Z^{'\lambda_2}$ containing $x$ as described
above, but if at time $t$ we find that a ball of radius $R$ is
contained in $Z^{'\lambda_2}[L_t^{'\lambda_2}(x)]$ we stop the
process. Let $T$ denote the random time at which the process stops.
Note that $T<\infty$ a.s., since if $Z^{'\lambda_2}(x)$ is unbounded,
then $Z^{'\lambda_2}(x)$ contains balls of radius $R$ a.s. by Lemma
\ref{ballemma}.  Let $F_1$ be the event that the process stops when a
ball of radius $R$ is found, and note that $Z^{'\lambda_2}(x)$ is a.s.
bounded on $F_1^c$. On $F_1$, we may (in some way independent of
$Z^{''}$) pick a point $Y$ such that $S(Y,R)$ is covered by
$Z^{'\lambda_2}[L_T^{'\lambda_2}(x)]$.

For $i=1,2$ let \[X^{\lambda_i}:=(X^{'\lambda_i}\cap
L_{T}^{'\lambda_2}(x))\cup(X^{''\lambda_i}\cap
L_T^{'\lambda_2}(x)^c)\] and $Z^{\lambda_i}:=\cup_{x\in
X^{\lambda_i}}S(x,1)$. In this way, $Z^{\lambda_i}$ is a Poisson
Boolean model with intensity $\lambda_i$ for $i=1,2$, and any ball
present in $Z^{\lambda_1}$ is also present in $Z^{\lambda_2}$.

Now put \[F_2:=F_1\cap \{S(Y,R)\mbox{ intersects an unbounded
  component of }Z^{''\lambda_1}\}.\]

But on $F_2$ there is some point in $Z^{\lambda_2}(x)$ which is at
distance less than or equal to two from some unbounded
$Z^{\lambda_1}$ component, that is $\{\tilde{D}(x)\le2\}$ occurs for
$Z$ so that $E_2(x)$ does not occur for $Z$. Since $E_2(x)$ is up to
a set of measure $0$ contained in $F_1$ we have that
\[\P[E_2(x)]\le\P[F_1\cap F_2^c].\] Since $Z^{'}$ and $Z^{''}$ are
independent it follows that
\[\P[F_2|F_1]=\P[S(Y,R)\mbox{ intersects an unbounded component of
}Z^{''\lambda_1}]\ge 1-a\] and consequently \[\P[F_1\cap
F_2^c]\le\P[F_2^c|F_1]<a.\] \noindent Since we may choose $a$
arbitrary small it follows that $\P[E_2(x)]=0$ as desired.

Next we argue that $\P[E_2(x)]=0$ for all $x$ implies $\P[E_2]=0$.
Let $D$ be a countable dense subset of $M$. Then $\P\left[\cup_{x\in
D}E_2(x)\right]=0$. But if $E_2$ occurs then $E_2(x)$ occurs for all
$x$ in some unbounded component of $C^{\lambda_2}$, in particular
for some $x$ in $D$, so it follows that $\P[E_2(x)]=0$ implies
$\P[E_2]=0$.

Next we show that $\P[E_1(x)]=0$. Let $E_1^f(x)$ be the event that
$E_1(x)$ occurs and all points in the $\lambda_1$-unstable unbounded
$C^{\lambda_2}$-component of $x$ which are at distance less than or
equal to three from some unbounded $C^{\lambda_1}$-component are
contained in the ball $S(0,N)$ for some random finite $N$. Let
$E_1^{\infty}(x)$ be the event that $E_1(x)$ occurs but that there
is no such finite $N$. Let $E_1^{f}$ and $E_1^{\infty}$ be the
events that $E_1^f(x)$ and $E_1^{\infty}(x)$ respectively happen for
some $x$.

First we show that $\P[E_1^f]=0$. Let $E_1^{f,M}:=E_1^f\cap\{N\le
M\}$. We will show that $\P[E_1^f]>0$ implies that $\P[E_2]>0$. So
suppose $\P[E_1^f]>0$. Then we may pick $M$ so large that
$\P[E_1^{f,M}]>0$. Again let $Z^{'}$ and $Z^{''}$ be independent
with the same distribution as $C$. Then for $i=1,2$ let
$Z^{\lambda_i}$ be the union of all balls from $Z^{'\lambda_i}$
centered within $S(0,M+1)$ together with the union of all balls from
$Z^{''\lambda_i}$ centered within $S(0,M+1)^c$. Then if
$\{Z^{'\lambda_2}[S(0,M+1)]=\emptyset\}$ occurs and $E_1^{f,M}$
occurs for $Z^{''}$ then $E_2$ occurs for $Z$. So since $Z^{'}$ and
$Z^{''}$ are independent we get
\[\P[E_2]\ge\P\left[Z^{'\lambda_2}[S(0,M+1)]=\emptyset\right]\P\left[E_1^{f,M}\right]>0\]
which is a contradiction, so $\P[E_1^f]=0$.

Finally we show that $\P[E_1^{\infty}]=0$. However the event
$E_1^{\infty}(x)$ is very similar to the event $B^{\infty}(x)$ in
the first proof of Theorem \ref{uniqmon1}, and is shown to have
probability $0$ in the same way. In the same way it then follows
that $\P[E_1^{\infty}]=0$. {$\Box$ \vspace{1em}}

\section{Connectivity}

In this section we show how $\lambda_u$ can be characterized by the
connectivity between big balls. This result will be used when we
study the model at $\lambda_u$ on a product space in the next
section. Let
\[\lambda_{BB}:=\inf\{\lambda:\,\lim_{R\rightarrow
  \infty}\inf_{x,y}\P[S(x,R)\leftrightarrow S(y,R)\mbox{ in }C^{\lambda}]=1\}.\] Note that
obviously $\lambda_{BB}\ge \lambda_c$. We will show the following:

\begin{thm}\label{ballconnect}
  For the Poisson Boolean model on a homogeneous space with $\lambda_u<\infty$ we have $\lambda_u=\lambda_{BB}$.
\end{thm}

The discrete counterpart of this result is Theorem 3.2 of
\cite{schonmann}, and the proof is similar. The proof is also
similar to the second proof of Theorem \ref{uniqmon1} above. First
we show that $\lambda_u\le \lambda_{BB}$.

\begin{proof}
Suppose that $\lambda_{BB}<\lambda_1<\lambda_2$. We will show that
at $\lambda_2$ there is a.s. a unique unbounded component. For
$i=1,2$ let
\[A_i(x,y):=\{\mu(C^{\lambda_1}(x))=\infty,\,\mu(C^{\lambda_1}(y))=\infty,\,C^{\lambda_i}(x)\neq
C^{\lambda_i}(y)\},\] and let \[A_i:=\bigcup_{x,y}A_i(x,y).\]

Since $\lambda_{BB}\ge \lambda_c$ we have by Theorem \ref{uniqmon1}
that any unbounded $\lambda_2$ component a.s. intersects some
unbounded $\lambda_1$ component. Therefore
\begin{equation}\label{inklusion1}
 \bigcup_{x,y} \{\mu(C^{\lambda_2}(x))=\infty,\,\mu(C^{\lambda_2}(y))=\infty,\,C^{\lambda_2}(x)\neq
C^{\lambda_2}(y)\}\subset A_2\cup N
\end{equation}
where $N$ is a set of measure $0$. In the same way as in the
second proof of Theorem \ref{uniqmon1} we have that
$\P[A_i(x,y)]=0$ for all $x$ and $y$ implies $\P[A_i]=0$. By
\ref{inklusion1}, $\P[A_2]=0$ implies $\P[\mbox{there is a unique
unbounded component at level }\lambda_2]=1$. Hence it is enough to
show that $\P[A_2(x,y)]=0$ for all $x$ and $y$.
\begin{definition}
  Suppose $C_1$ and $C_2$ are two distinct components in the
  Poisson Boolean model. A pair of Poisson points $x_1\in C_1$ and
  $x_2\in C_2$ is called a \emph{boundary-connection} between $C_1$
  and $C_2$ if $d(x_1,x_2)<4$ (so that the distance between their
  corresponding balls is $<2$) or there is a sequence of
  Poisson-points $y_1,\dots, y_n$ such that
\begin{itemize}
\item the ball centered around $y_i$ intersects the ball centered around $y_{i+1}$
for all $i$.
\item $y_i$ is outside $C_1$ and $C_2$ for all $i$.
\item $d(x_1,y_1)<4$ and $d(x_2,y_n)<4$.
\end{itemize}
\end{definition}

Note that if there is a boundary connection between two components,
then at most two more balls are needed to merge them into one
component.

If $x,y\in C^{\lambda_1}$ and $C^{\lambda_1}(x)\neq C^{\lambda_1}(y)$,
let $B(x,y)$ be the number of boundary connections between
$C^{\lambda_1}(x)$ and $C^{\lambda_1}(y)$. Let
\[A_1^0(x,y):=A_1(x,y)\cap \{B(x,y)=0\},\] \[A_1^f(x,y):=A_1(x,y)\cap
\{B(x,y)<\infty\},\] \[A_1^\infty(x,y):=A_1(x,y)\cap
\{B(x,y)=\infty\},\] and for $t\in\{0,f,\infty\}$ let $A_1^t$ be the
event that $A_1^t(x,y)$ happens for some $x$ and $y$. In the same way
as before it is seen that $\P[A_1^t(x,y)]=0$ for all $x$ and $y$
implies $\P[A_1^t]=0$.

Next we will argue that
\begin{equation}
\label{bondconnect1} \P[A_1^0(x,y)]=0\mbox{ for all }x\mbox{ and }y.
\end{equation}
Let $Z^{'\lambda_1}$ and $Z^{''\lambda_1}$ be two independent copies
of the Poisson Boolean model with intensity $\lambda_1$ and let
$X^{'\lambda_1}$ and $X^{''\lambda_1}$ be their underlying Poisson
processes. Since $\lambda_1>\lambda_{BB}$ we may for any $a>0$ pick
$R=R(a)$ such that
\[\inf_{z_1,z_2}\P_{\lambda_1}[S(z_1,R)\leftrightarrow S(z_2,R)]>1-a.\]
Fix $x$ and $y$ and grow the component of $x$ in $Z^{'\lambda_1}$ (as
described earlier) but stop if a ball of radius $R$ is found. Do the
same for $y$. Let $F_1$ be the event that the processes are stopped
when balls of radius $R$ are found, and note that $A_1^0(x,y)$ is up
to a set of measure $0$ included in $F_1$. Let $T_x$ and $T_y$ denote
the random times at which the processes are stopped.  On $F_1$ we pick
$X$ and $Y$ in some way independent of $Z^{''\lambda_1}$ such that
$S(X,R)\subset Z^{'\lambda_1}[L_{T_x}^{'\lambda_1}(x)]$ and
$S(Y,R)\subset Z^{'\lambda_1}[L_{T_y}^{'\lambda_1}(y)]$. Let
\[X^{\lambda_1}:=(X^{'\lambda_1}\cap(L_{T_x}^{'\lambda_1}(x)\cup L_{T_y}^{'\lambda_1}(y)))\cup
(X^{''\lambda_1}\cap (L_{T_x}^{'\lambda_1}(x)\cup
L_{T_y}^{'\lambda_1}(y))^c)\] and $Z^{\lambda_1}:=\cup_{x\in
X^{\lambda_1}}S(x,1)$. The distribution of $Z^{\lambda_1}$ is by
construction the distribution of the Poisson Boolean model with
intensity $\lambda_1$. Put
\[F_2:=F_1\cap \{S(X,R)\leftrightarrow S(Y,R) \mbox{ in
}Z^{''\lambda_1}\}.\] If we are on $F_2$ then either
$\{Z^{\lambda_1}(x)=Z^{\lambda_1}(y)\}$ occurs or $\{B(x,y)\ge 1\}$
occurs and in neither case we are on $A_1^0(x,y)$. Since
\[\P[F_2|F_1]=\P[S(X,R)\leftrightarrow S(Y,R)\mbox{ in
}Z^{''\lambda_1}]>1-a\] it therefore follows that
\[\P[A_1^0(x,y)]\le\P[F_1\cap F_2^c]\le \P[F_2^c|F_1]<a\] proving
(\ref{bondconnect1}).

Next we show that
\begin{equation}
\label{bondconnect2} \P[A_1^f]=0.
\end{equation}

Let $A_1^{f,N}$ be the event there are two distinct unbounded
components in $C^{\lambda_1}$ such there are a finite number of
boundary connections between them and they are all contained in the
ball $S(0,N)$ for some random finite $N$. Suppose $\P[A_1^f]>0$ and
pick $N$ so large that $\P[A_1^{f,N}]>0$. Let $Z^{\lambda_1}$ be the
union of the balls from $Z^{'\lambda_1}$ centered outside $S(0,N)$
and the balls from $Z^{''\lambda_1}$ centered inside $S(0,N)$. Now
suppose that $A_1^{f,N}$ happens for $Z^{'\lambda_1}$ and that
$\{Z^{''\lambda_1}[S(0,N)]=\emptyset\}$ happens. Then we can find
two points $\tilde{x}$ and $\tilde{y}$ in separate unbounded
components of $Z^{\lambda_1}$ such that there are no boundary
connections between them. It follows by the independence of $Z^{'}$
and $Z^{''}$ that
\[\P[A_1^0]\ge
\P\left[A_1^{f,N}\right]\P\left[Z^{''\lambda_1}[S(0,N)]=\emptyset\right]>0,\]
a contradiction which proves (\ref{bondconnect2}).

Now if $A_1^{\infty}(x,y)$ happens, then there are infinitely many
boundary connections between $C^{\lambda_1}(x)$ and
$C^{\lambda_1}(y)$ and a.s. no bounded region contains all boundary
connections. Therefore $C^{\lambda_1}(x)$ and $C^{\lambda_1}(y)$
will almost surely have been merged into one unbounded component at
level $\lambda_2$ by balls that appear in the coupling between level
$\lambda_1$ and $\lambda_2$. So $\P[A_2(x,y)|A_1^{\infty}(x,y)]=0$.
Thus, since $A_2(x,y)\subset A_1(x,y)$ and $A_1(x,y)$ is partitioned
by $A_1^f(x,y)$ and $A_1^{\infty}(x,y)$ we conclude
\[\P[A_2(x,y)]=\P[A_2(x,y)|A_1^f(x,y)]\P[A_1^f(x,y)]+\P[A_2(x,y)|A_1^{\infty}(x,y)]\P[A_1^{\infty}(x,y)]=0,\]
for all $x$ and $y$ and so $\lambda_u\le \lambda_{BB}$.

Next we show the easier result that $\lambda_u\ge \lambda_{BB}$.
Suppose $\lambda>\lambda_u$. By Theorem \ref{uniqmon1} there is a.s.
a unique unbounded component in $C^{\lambda}$ which we denote by
$C_{\infty}^{\lambda}$. By the continuum version of the FKG
inequality (see \cite{meester}) and the fact that there is an
isometry mapping $x$ to $y$ it follows that
\[\begin{split}\P_{\lambda}[S(x,R)\leftrightarrow
S(y,R)]& \ge\P_{\lambda}[S(x,R)\mbox{ and }S(y,R)\mbox{ intersects
}C_{\infty}^{\lambda}]\\& \ge\P_{\lambda}[S(x,R)\mbox{ intersects
}C_{\infty}^{\lambda}]^2.\end{split}\] Since $\lim_{R\rightarrow
\infty}\P_{\lambda}[S(x,R)\mbox{ intersects
}C_{\infty}^{\lambda}]=1$ it follows that $\lambda>\lambda_{BB}$ and
thus $\lambda_u\ge \lambda_{BB}$.
\end{proof}

\section{The situation at $\lambda_u$ on $\hyp \times \R$}

This section is devoted to the proof of Theorem \ref{htimesrsats}.  We
introduce some new notation: if the points $x,y\in \hyp\times \R$ are
in the same component of $C^{\lambda}$ then $d_{X^{\lambda}}(x,y)$ is
the smallest number of balls in that component forming a sequence that
connects $x$ to $y$. For a set $A$ we let $C^{\lambda}(A)$ be the
union of all components of $C^{\lambda}$ that intersect $A$. The
length of a curve $\gamma\subset\hyp$ will be denoted by $L(\gamma)$.
In this proof $\mu=\mu_{\hyp}$ and $d=d_{\hyp\times {\mathbb R}}$.

\begin{proof}
  As noted earlier, it is the case that $\lambda_u(\hyp\times
  \R)<\infty$. Suppose that $\lambda_*$ is such that there is a.s. a
  unique unbounded component in the Poisson Boolean model with
  intensity $\lambda_*$ on $\hyp\times \R$. We consider the monotone
  coupling of the model for all intensities below $\lambda_*$. We will
  show that there is some intensity below $\lambda_*$ that also a.s.
  produces a unique unbounded component. Denote the unbounded
  component at $\lambda_*$ with $C^{\lambda_*}_{\infty}$. For any
  $r>0$, any positive integer $n$, and any $\lambda\in(0,\lambda_*)$
  we define the following three random sets:
\[A_1(r):=\{z\in\hyp\times \R\,:\,S(z,r)\cap
C_{\infty}^{\lambda_*}\neq\emptyset\}\]\[A_2(r,n):=\{z\in\hyp\times
\R\,:\,\sup\{d_{X^{\lambda_*}}(s,t)\,:\,s,t\in S(z,r+1/2)\cap
C_{\infty}^{\lambda_*}\}<n\}\]\[A_3(r,n,\lambda):=\{z\in\hyp\times
\R\,:\,S(z,r+2n)\cap(X^{\lambda_*}\setminus
X^{\lambda})=\emptyset\}.\] Then put
\[A(r,n,\lambda):=A_1(r)\cap A_2(r,n)\cap
A_3(r,n,\lambda).\] Pick $y_1,\,y_2\in {\mathbb R}$ and let
\[D:=D(y_1,y_2,r,n,\lambda)=\{x\in \hyp\,:\,(x,y_1)\in
A(r,n,\lambda)\mbox{ and } (x,y_2)\in
A(r,n,\lambda)\}.\] Then $D$ is a random set in $\hyp$ such
that the law of $D$ is Isom($\hyp$)-invariant. Next we will show
that we can choose the parameters $r,n$ and $\lambda$ in such a way
that $D$ contains unbounded components with positive probability.

To do this, we let $\tilde{C}$ be a Poisson Boolean model in $\hyp$
with intensity $\tilde{\lambda}$. Let $B$ be a bounded connected set in
$\hyp$. Choose $\tilde{\lambda}$ so big that $\hyp$, we have
$\E[L(B\cap \partial \tilde{C})]<\E[\mu(B\cap \tilde{C})]$. By Lemma
5.2 in \cite{tykesson}, $\tilde{C}$ contains unbounded components with
probability $1$. Let $\tilde{C}^D=\tilde{C}^D(y_1,y_2,r,n,\lambda)$ be
the union of all balls in $\tilde{C}$ that are completely covered by
$D$.

Suppose $E$ is some bounded connected set in $\hyp\times{\mathbb R}$.
It is clear that
\begin{equation}\label{grans1}\lim_{r\rightarrow \infty}\P[E\subset
  A_1(r)]=1,\end{equation} and that for fixed $r$,
\begin{equation}\label{grans2}\lim_{n\rightarrow \infty}\P[E\subset
  A_2(r,n)]=1,\end{equation} and that for fixed $r$ and $n$,
\begin{equation}\label{grans3}\lim_{\lambda\uparrow\lambda_0}\P[E\subset
  A_3(r,n,\lambda)]=1.\end{equation}

Put $\delta:=\E[\mu(B\cap \tilde{C})]-\E[L(B\cap \tilde{C})]$. By
(\ref{grans1}), (\ref{grans2}) and (\ref{grans3}) we get that we can
find first $r_1$ big enough, and then $n_1$ big enough, and finally
$\lambda_1$ close enough to $\lambda_*$ so that $\E[\mu(B\cap
\tilde{C})]-\E[\mu(B\cap \tilde{C}^D)]<\delta/2$ and $\E[L(B\cap
\partial \tilde{C}^D)]-\E[L(B\cap \partial \tilde{C})]<\delta/2.$ With
these choices of parameters, $\E[\mu(B\cap \tilde{C}^D)]>\E[L(B\cap
\tilde{C}^D)]$, so by Lemma 5.2 in \cite{tykesson}, we get that
$\tilde{C}^D$ contains unbounded components with positive probability.
Since $\tilde{C}^D\subset D$, this implies that $D$ contains unbounded
components with positive probability. Since the event that $D$
contains unbounded components is Isom($\hyp$)-invariant and determined
by the underlying Poisson processes in the model, $D$ contains
unbounded components with probability $1$.

So we can find an infinite sequence of points $u_1,u_2,...\in \hyp$
such that they are all in the same component of $D$,
$d(u_i,u_{i+1})<1/2$ for all $i$ and $d(u_1,u_i)\rightarrow \infty$ as
$i\rightarrow \infty$.  Since $(u_i,y_1)\in A_1$ there is some ball
$s_i$ in $C^{\lambda_0}_{\infty}$ centered within distance $r_1+1$
from $(u_i,y_1)$. Since $d((u_i,y_1),(u_{i+1},y_1))<1/2$ and
$(u_i,y_1)\in A_2$ for all $i$ there is a sequence of at most $n$
balls in $C^{\lambda_0}_{\infty}$ connecting $s_i$ to $s_{i+1}$. Since
the distance between the center of any ball in this sequence and
$(u_i,y_1)$ is at most $r_1+2n$ and $(u_i,y_1)\in A_3$, all balls in
the sequence is present also at level $\lambda_1$. Thus there is an
unbounded component in $C^{\lambda_1}$ that comes within distance
$r_1$ from $(u_i,y_1)$ for all $i$. In the same way there is an
unbounded component in $C^{\lambda_1}$ that comes within distance
$r_1$ from $(u_i,y_2)$ for all $i$.

 Now choose $\lambda_2$ and $\lambda_3$ so that
 $\lambda_1<\lambda_2<\lambda_3<\lambda_*$. For $x\in\hyp$ let
 $D(x)$ be the component of $D$ containing $x$. Then we have from the
 above that
\begin{equation}\label{ball1}\P[S((x,y_1),r_1)\leftrightarrow
  S((x,y_2),r_1)\mbox{ in }C^{\lambda_2}|\mu(D(x))=\infty]=1.\end{equation} This follows from
the fact that the two unbounded components at level $\lambda_1$ above will almost surely
be connected by balls appearing in the coupling between level
$\lambda_1$ and $\lambda_2$. Fix $a$ small and let $r_2$ be such
that for $x\in\hyp$ the ball $S(x,r_2)$ in $\hyp$ intersects an unbounded
component of $D$ with probability at least $1-a/2$. Let $R=r_1+r_2$.
If $S(x,r_2)$ intersects an unbounded component of $D$ then by
(\ref{ball1}) it follows that a.s.
$S((\tilde{x},y_1),r_1)\leftrightarrow S((\tilde{x},y_2),r_1)$ in
$C^{\lambda_2}$ for some for some point $\tilde{x}\in \hyp$ such
that $d_{\hyp}(x,\tilde{x})\le r_2$, so $S((x,y_1),R)\leftrightarrow
S((x,y_2),R)$ in $C^{\lambda_2}$. Thus
\begin{equation}\label{ball2}\P[S((x,y_1),R)\leftrightarrow
  S((x,y_2),R)\mbox{ in }C^{\lambda_2}]\ge 1-a/2.\end{equation} Fix
two points $z_1=(u_1,v_1)$ and $z_2=(u_2,v_2)$ of $\hyp\times \R$.
For $y\in \R$ let \[F_y:=\{S(z_1,R)\leftrightarrow
S((u_1,y),R)\mbox{ in }C^{\lambda_2}\}\cap \{S(z_2,R)\leftrightarrow
S((u_2,y),R)\mbox{ in }C^{\lambda_2}\}\] By (\ref{ball2}) we get
$\P[F_y]\ge 1-a$ for all $y$. In particular it follows that with
probability at least $1-a$ the set $\{y\in {\mathbb R}\,:\:F_y\mbox{
  occurs }\}$ is unbounded. But then the set of points in
$C^{\lambda_2}(S(z_1,R))$ that come within distance
$2R+d_{\hyp}(u_1,u_2)$ from $C^{\lambda_2}(S(z_2,R))$ is unbounded.
But if this occurs then some component in $C^{\lambda_2}$ intersecting $S(z_1,R)$ will a.s. be connected to some component in $C^{\lambda_2}$ intersecting $S(z_2,R)$ by balls occurring in the coupling between level
$\lambda_2$ and $\lambda_3$. That is,
\[\P[S(z_1,R)\leftrightarrow S(z_2,R)\mbox{ in }C^{\lambda_3}]\ge 1-a.\]
Since $a$ is arbitrary small it follows by Theorem \ref{ballconnect}
there is a.s. a unique unbounded component in $C^{\lambda_3}$.
\end{proof}

\noindent {\bf Remark.} Of course, there is nothing special about
${\mathbb R}$ in the proof of Theorem \ref{htimesrsats}. The proof
works without any modifications if ${\mathbb R}$ is replaced by any
noncompact homogeneous space $M$ such that $\lambda_u(\hyp\times M)<\infty$.
Also, it is possible to show a version of Lemma 5.2 in \cite{tykesson}
for ${\mathbb H}^n$ for any $n\ge 3$. Therefore Theorem \ref{htimesrsats}
holds for $\hypn\times M$ for any $n\ge 2$ and any noncompact homogeneous space
if $\lambda_u(\hypn\times M)<\infty$.

\section{Further problems}

In this section we list some open problems.  

1. For which manifolds is $\lambda_u<\infty$?

2. In \cite{tykesson} it is shown that
$\lambda_c(\hypn)<\lambda_u(\hypn)$ for any $n\ge 2$ if the radius of
the percolating balls is big enough (for $n=2$ this is shown for any
radius). For which manifolds is $\lambda_c<\lambda_u$?

3. For which manifolds with $\lambda_u<\infty$ is there a.s. a unique
unbounded component at $\lambda_u$? For which manifolds is there a.s.
not a unique unbounded component at $\lambda_u$?

\bigskip\noindent {\bf Acknowledgement:} I want to thank Johan
Jonasson, my advisor, for useful discussions and comments.

\end{document}